\documentclass[journal]{IEEEtran}

\ifCLASSINFOpdf
   \usepackage[pdftex]{graphicx}
\else
   \usepackage[dvips]{graphicx}
\fi

\usepackage[cmex10]{amsmath}
\usepackage{amsfonts}
\usepackage{amssymb}

\usepackage{mathtools}
\usepackage[noend,linesnumbered]{algorithm2e}
\usepackage{makecell}

\usepackage{booktabs} 
 \usepackage{csvsimple}
\usepackage{multirow}
\usepackage{siunitx}
\usepackage{url}
\newrobustcmd{\prevscenarios}{}

\usepackage{xcolor}
\newcommand{\mv}[1]{\boldsymbol{#1}}

\usepackage[normalem]{ulem}
\newcommand{\stkout}[1]{\ifmmode\text{\sout{\ensuremath{#1}}}\else\sout{#1}\fi}

\usepackage{soul}    

\begin{document}
\title{Machine Learning-Enabled Large-Scale Capacity Expansion Planning under Uncertainty}

\author{\IEEEauthorblockN{Taehyeon Kwon\IEEEauthorrefmark{1},
Anirudh Subramanyam\IEEEauthorrefmark{1}}
\IEEEauthorblockA{\IEEEauthorrefmark{1}Department of Industrial and Manufacturing Engineering\\
The Pennsylvania State University,
University Park, PA 16802, USA}
\thanks{\noindent This material is based upon work supported by the United States Department of Energy, award number DE-SC0023361.}

}

\maketitle

\begin{abstract}
Capacity expansion planning under uncertainty requires selecting a scenario count and representative operational horizon to estimate average production costs. Small choices risk unreliable plans, while large choices become intractable. %
We propose AutoSCEP, an automated, statistically grounded procedure that, for a fixed plan, selects the minimum sufficient scenario count and horizon length to estimate production costs to a given precision. Using these estimates, we train linear and neural surrogates to approximate expected production costs for arbitrary plans, and embed the surrogates within the planning model. On the continental-scale EMPIRE system, AutoSCEP attains 2\% optimality gap on a reduced model and 8\% gap on a large model, outperforming parallel progressive hedging under equal wall-clock budgets that include data generation, training, and solve times. Where the reduced model's optimum is available, investment patterns broadly align with the benchmark. Our approach enables high-resolution uncertainty modeling at realistic system scales.
\end{abstract}

\begin{IEEEkeywords}
capacity expansion planning, automated parameter selection, surrogate modeling, machine learning
\end{IEEEkeywords}

\section{Introduction}

Generation and transmission capacity planning under uncertainty is one of the most challenging problems in modern power systems. This challenge stems from multiple sources: rapidly growing load driven by data centers, %
aging infrastructure,
increasing uncertainty from the integration of renewable energy sources to meet policy mandates such as carbon emission reduction targets, and emerging new technologies \cite{lumbreras2016new, cao2025review}. These challenges can be addressed by solving Capacity Expansion Planning (CEP) problems,
which co-optimize capital and production costs to determine the least-cost future capacity expansion plan while satisfying physical, technical, and policy constraints across numerous uncertainty realizations \cite{conejo2016investment}.

Stochastic programming (SP) is one of the most widely used techniques to model and solve CEP problems under uncertainty,  %
since it effectively captures uncertainties by modeling them as scenarios that can represent either historical data, or samples from a probability distribution, or expert-driven narratives \cite{roald2023power}. However, their computational complexity rapidly grow with respect to model resolution. %
Low-resolution models may produce unreliable or costly investment plans, whereas high-resolution models may yield more reliable plans but are typically computationally intractable \cite{roh2009market,poncelet2016selecting}. 
This dilemma is further complicated by the need to explicitly account for uncertainties in the planning process and the multiscale nature of CEP problems that need to determine long-term investment plans (over years to decades) while considering future short-term operational variability to estimate production costs \cite{krishnan2016co}.

In SP-based formulations of CEP problems, model resolution is dictated by the number of uncertain scenarios and the operational horizon. The number of scenarios is a measure of how many possible future outcomes are considered when making long-term investment plans. The operational horizon refers to the resolution of each representative operational period (e.g., days, hours, or minutes) and the number of such chronological or consecutive periods. It provides a measure of operational variability, with longer and finer horizons providing more accurate estimates of future production costs \cite{bylling2020impact}. 
Software (e.g., PLEXOS) commonly used by transmission system operators (TSOs) often simply use single hourly slices to maintain tractability.

There are two main approaches for handling large scenario counts and operational horizons. The first approach is scenario reduction and clustering \cite{dupavcova2003scenario}. 
Existing methods sample from large Monte Carlo or historical datasets and then compress them to a small number for tractability purposes. Examples include random or historical sampling of a single uncertainty \cite{feng2013scenario}, scenario reduction by clustering \cite{defourney2013}, and multi-scale selection of strategic and operational states \cite{liu2017hierarchical}. Clustering metrics span raw multivariate input patterns %
and power flow or congestion patterns \cite{alvarez2017,ploussard2017,sun2018}.

The second approach employs decomposition algorithm to solve SP formulations involving a large number of scenarios. Benders decomposition %
has been extended for capacity expansion planning problems \cite{zhang2024decomposition}. However, for large-scale networks with long planning horizons, it can require substantial computational effort \cite{jacobson2024computationally}. Progressive hedging (PH) has also been widely used for capacity planning problems 
\cite{munoz2015scalable, kaisermayer2021progressive}. It decomposes subproblems by scenario and solves them in parallel using an augmented Lagrangian method to enforce non-anticipativity constraints. %

Recently, a new solution approach has emerged to tackle scalability and computational issues. Neural Two-Stage Stochastic Programming (Neur2SP) \cite{patel2022neur2sp} uses a neural network surrogate model with ReLU activation layers to approximate the expected second-stage value function, which can be embedded as mixed integer linear constraints \cite{fischetti2018deep}. Unlike other solution algorithms, Neur2SP is a `one-shot' approach; after initial training, the neural network surrogate remains fixed and can be reused across multiple
instances with varying first-stage parameters without retraining~\cite{patel2022neur2sp}. This offers a potential path to address the high temporal and uncertainty resolutions in CEP problems. 
The method has been extended and improved along multiple directions \cite{kronqvist2023alternating,alcantara2025quantile,liu2025icnn}.

However, despite their potential, adapting these methods to CEP problems under uncertainty is not straightforward. For example, generating data for training the machine learning surrogates would entail solving millions of second-stage operational production costing problems. This can easily take several days, since estimating the expected production cost requires large scenarios counts over long operational horizons.
Existing methods require that these parameters be fixed \textit{a priori} and any changes to their values would entail either retraining the machine learning model or re-solving traditional algorithms (such as PH) from scratch, which can become highly time and resource intensive.

We propose AutoSCEP, short for Automated Surrogate-based Capacity Expansion Planning, a Neur2SP-inspired framework that generates  capacity investment plans without requiring any manual selection of scenario counts and operational horizon lengths. It only requires access to a scenario generator that can generate operational uncertain scenarios for any given horizon length. Our specific contributions are:
\begin{itemize}
    \item
    We develop an adaptive parameter selection algorithm, which, for a fixed capacity investment plan, uses statistical significance testing to automatically select the minimum scenario count and representative operational horizon length needed to estimate average production costs to a prescribed precision. The algorithm is self-contained and can be used independently to evaluate the expected cost of candidate investment plans. %

    \item We develop a complete training pipeline that combines the adaptive parameter selection algorithm %
    with constraint propagation-based sampling of feasible capacity expansion plans
    to efficiently generate training data. Using these data, we train linear and neural network surrogates that approximate expected production costs for arbitrary expansion plans, and then embed these surrogates within the planning model to obtain an optimal investment plan.
    
    \item We conduct extensive experiments on the European EMPIRE system \cite{backe2022empire}, a large-scale planning model over 50 years. %
    We evaluate AutoSCEP by comparing it against (i) traditional SP methods, including the classical extensive form, Benders decomposition, and progressive hedging, and (ii) linear and neural network surrogates that use \textit{a priori} fixed parameters during training.
    We find that under equal wall-clock budgets that include data generation, training, and solve times, AutoSCEP can outperform parallel progressive hedging.
\end{itemize}
All code and data used in this study are publicly available at \url{https://github.com/Subramanyam-Lab/AutoSCEP}.

This paper is organized as follows: Section~\ref{sec:formulation} details the planning model we consider. Section~\ref{sec:methodology} presents our proposed solution methodology, including the data generation and ML surrogate training and embedding. Section~\ref{sec:experiments} describes the experimental setup. Section~\ref{sec:results} analyzes comparative results, and Section~\ref{sec:conclusions} concludes with a summary of our findings.

\section{Stochastic Capacity Expansion Planning Problem Formulation} \label{sec:formulation}

We consider the following two-stage stochastic programming formulation of the capacity expansion problem:
\begin{equation}
\label{eq:2sp_cep}
\begin{aligned}
\mathop{\mathrm{minimize}} \;\quad & \mv{c}^{\top}_{\mathrm{inv}}\mv{x} + \mathbb{E}_{\mv{\omega}} [Q(\mv{x}, \mv{\omega})] \\
\mathrm{subject\ to} \quad & \mv{A}\mv{x} \leq \mv{b}, \; \mv{\ell} \leq \mv{x} \leq \mv{u}, 
\end{aligned}
\end{equation}
where
$\mv{x} = (x^{\theta}_{i} )_{i \in S_\theta, \theta \in \Theta}$
collects capacity investments (MW) in various technology types (e.g., $\Theta = \{\text{gen}, \text{tran}, \text{stor}\}$), and $S_\theta$ consists of candidate investments in technology type $\theta \in \Theta$.
For example, $S_{\text{gen}}$ may consist of `coal in zone \emph{A}' or `solar in zone~\emph{B}', $S_{\text{tran}}$ may consist of `HVDC line from zone \emph{C} to zone~\emph{D}', and $S_{\text{stor}}$ may consist of `pumped hydro in zone~\emph{E}'. %

The second-stage cost function $Q(\mv{x}, \mv{\omega})$ captures production costs under investment decisions $\mv{x}$ and future scenario $\mv{\omega} \in \Omega$:
\begin{equation}
	\label{eq:second_stage}
	Q(\mv{x}, \mv{\omega}) = \min_{\mv{y} \in \mathcal{Y}(\mv{x}, \mv{\omega})} \sum_{h\in \mathcal{H}} c_{\text{prod},\mv{\omega}}(\mv{y}_{h\mv{\omega}}),
\end{equation}
where $\mathcal{Y}(\mv{x}, \mv{\omega})$ denotes the set of feasible operational decisions $\mv{y}_{h\mv{\omega}}$ (e.g., hourly generation levels, transmission line flows etc.), defined by power balance equations, generation capacity constraints, and ramping limits (besides others), over an operational time horizon $\mathcal{H}$.
The latter is a set of chronological hours over which power grid operations are modeled, which can vary from hourly slices $|\mathcal{H}| = 1$ to annual production cost models with $|\mathcal{H}| = 8760$.
In particular, $|\mathcal{H}| > 1$ enables modeling of operational variability due to ramping limits and storage facilities \cite{bylling2020impact,poncelet2016selecting}.
The overall objective function in~\eqref{eq:2sp_cep} minimizes the sum of annualized investment costs and (present worth of) future expected production costs.

Several remarks about the model are in place.
First, the generation, transmission, and storage capacity investments $\mv{x}$ are modeled as continuous variables (where any nonnegative MW value can be chosen for an entire class of technologies at a given location), rather than binary variables for individual facility buildouts.
This aligns with existing large-scale models used for policy analyses, such as TIMES \cite{loulou2005documentation}, ReEDS \cite{ho2021reeds}, and EMPIRE \cite{backe2022empire}, as well as industry software (e.g., PLEXOS) that use zonal network representations.

Second, the operational model $Q(\mv{x}, \mv{\omega})$ can be quite flexible as long as it satisfies complete recourse, i.e., is feasible for any $\mv{x}$ and $\mv{\omega}$, which can be ensured by incorporating load shedding decisions.
Apart from that, it may use either transshipment representations of transmission \cite{backe2022empire} or linearized DC load flows of aggregations of transmission facilities \cite{ho2021reeds} or even binary unit commitment decisions \cite{garcia2023computational}.
Similarly, the uncertainties $\mv{\omega}$ may include both short-term fluctuations, such as generator availability, and load/wind/solar profiles, as well as long-term load growth and fuel prices, besides others. 

Finally, as written, \eqref{eq:2sp_cep} models a single investment period for a target year in the future~\cite{conejo2016investment, boffino2019two}. In practice, capacity expansion requires determining when to invest, with investment decisions at yearly or five-year intervals over a multi-decadal horizon $\mathcal{T}$.
Although traditional multi-stage models become intractable in such cases,  multi-horizon stochastic programming \cite{kaut2014multi} gives a tractable approximation:
\begin{equation}
	\label{eq:mhsp_cep}
	\begin{aligned}
		\min_{\mv{x}} \quad & \sum_{t \in \mathcal{T}} \mv{c}^{\top}_{\mathrm{inv},t}\mv{x}_t +  \sum_{t \in \mathcal{T}}  \mathbb{E}_{\mv{\omega}_t} [Q(\mv{v}_t, \mv{\omega}_t)]
	\end{aligned}
\end{equation}
Here, $x_{it}^\theta$ represents new capacity to be made available in investment period $t \in \mathcal{T}$, $v_{it}^{\theta} = \bar{x}^{\theta}_{it} + \sum_{t'=t-l_i}^t x_{it'}^{\theta}$ and $\bar{x}^{\theta}_{it}$ are the cumulative and pre-existing capacity available in period~$t$,  respectively, $l_i$ is the lifetime of infrastructure $i \in S_{\theta}$, and $Q$ has the same definition as~\eqref{eq:second_stage} but with additional cost discounting.
Crucially, \eqref{eq:mhsp_cep} retains the two-stage structure of~\eqref{eq:2sp_cep}, where all investment decisions $\{\mv{x}_t\}_{t \in \mathcal{T}}$ constitute the first stage, while operational decisions across all periods and scenarios constitute the second stage. Our proposed method applies to any model with this two-stage structure; in particular, all experimental results in this paper are based on an instantiation of the multi-horizon stochastic model~\eqref{eq:mhsp_cep}.

\section{ML Surrogate-Based Solution Methodology}\label{sec:methodology}

Evaluating the expected value function $\mathbb{E}_{\mv{\omega}} [Q(\mv{x}, \mv{\omega})]$ for fixed $\mv{x}$ is already computationally intractable~\cite{hanasusanto2016comment}.
Prior works~\cite{patel2022neur2sp,kronqvist2023alternating,liu2025icnn} use supervised learning to approximate the mapping $\left(\mv{x}, \{\mv{\omega}^s\}_{s=1}^S\right)  \mapsto \frac{1}{S} \sum_{s=1}^S Q(\mv{x}, \mv{\omega}^s)$ defined by the sample average approximation for some $S$ chosen \textit{a priori}.
For the large-scale capacity planning problems in our work, however, learning this mapping is extremely compute-intensive and practically untenable.
Therefore, we instead build a surrogate function $\Phi^{\mathrm{ML}}$ that directly approximates the expected value function, i.e., the mapping $\mv{x} \mapsto \mathcal{Q}(\mv{x}) \coloneqq\ \mathbb{E}_{\mv{\omega}} [Q(\mv{x}, \mv{\omega})]$ in~\eqref{eq:2sp_cep} or $\mv{x} \mapsto \sum_{t \in \mathcal{T}}\mathbb{E}_{\mv{\omega}_t} [Q(\mv{v}_t, \mv{\omega}_t)]$ in~\eqref{eq:mhsp_cep}.
Since $\mathcal{Q}(\mv{x})$ is a deterministic function of the first-stage decisions alone, its surrogate $\Phi^{\mathrm{ML}}(\mv{x})$ also does not depend on any particular scenario sample and can therefore be embedded in the planning model without a sample average approximation.

The methodology consists of the following steps.
\textit{(A)} Generate training data: $(\mv{x}^k, \widehat{\mathcal{Q}}^k)_{k=1}^K$, where $\mv{x}^k$ is a randomly sampled feasible first-stage investment decision and $\widehat{\mathcal{Q}}^k$ is its corresponding label, i.e., an estimator of the true expected production cost $\mathcal{Q}(\mv{x}^k)$ under $\mv{x}^k$. 
\textit{(B)} Define a surrogate model architecture and train it using the generated data.
\textit{(C)} Embed the trained model in the SP formulation~\eqref{eq:2sp_cep} to obtain the capacity expansion decisions.

\subsection{Generating Training Data}
Data generation involves two steps: \textit{(i)} sampling feasible first-stage decisions, and \textit{(ii)} generating training labels.  

It is important to generate a diverse set of feasible first-stage decisions so that the trained surrogate model will generalize across the entire decision space.
Prior works~\cite{patel2022neur2sp,kronqvist2023alternating,liu2025icnn} use ad-hoc problem-specific procedures based on independently sampling each variable within its bounds; however, doing so can cause the resulting capacity investment vector $\mv{x}$ to become infeasible with respect to coupling constraints (such as budgets) across periods $t \in \mathcal{T}$ and generation technologies.

We exploit ideas rooted in constraint propagation, also known as domain filtering in artificial intelligence \cite{davis1987constraint,hooker2007integrated}.
However, unlike these methods that were designed for optimization and converge only asymptotically, we instead use them for sampling.
Algorithm~\ref{alg:sequential_sampling} describes the procedure.
Once a variable value $x_i$ is sampled, the algorithm propagates the sampled value through every constraint in which it appears, tightening the bounds of all remaining (unsampled) variables that appear in a common constraint with $x_i$.
This ensures that the feasible region is dynamically narrowed so that any subsequent sampled variable is guaranteed to satisfy its constraints.
The algorithm converges in $O(Kmn)$ time, where $K$, $m$ and $n$ are the number of samples, constraints, and decision variables, respectively.
When $\mv{\ell}=\mv{0}$, and $\mv{A}$, $\mv{b}$, and $\mv{u}$ all have nonnegative entries, as is typical in CEP models, the algorithm also guarantees that the resulting samples are always feasible.

\RestyleAlgo{ruled}
\begin{algorithm}[hbt!]
\SetInd{0.35em}{0.5em}
\SetKwProg{Function}{function}{:}{end function}
\small
\caption{Constraint Propagation-Based Sampling}\label{alg:sequential_sampling}
\DontPrintSemicolon
\KwIn{sample size $K$, variable bounds $\mv{\ell}, \mv{u}\in\mathbb{R}^n$, constraint data $\mv{A} = (a_{rj}) \in \mathbb{R}^{m\times n}, \mv{b} \in \mathbb{R}^m$}
\KwOut{feasible sample set $X = \{\boldsymbol{x}^{1},\dots,\boldsymbol{x}^{K}\}$}
\Function{$f(r, i, \bar{J}, \mv{\ell}', \mv{u}', \mv{b}')$}{
\Return $\displaystyle \frac{b'_r}{a_{ri}} -
 \sum_{j\in \bar{J}: a_{rj} > 0}  \frac{a_{rj}\ell'_{j}}{a_{ri}}
-  \sum_{j\in \bar{J}: a_{rj} < 0}  \frac{a_{rj}u'_{j}}{a_{ri}}$\;
}
$X \gets \emptyset$\;
\For{$k = 1, \ldots, K$}{
$\boldsymbol{x}^{k} \gets \boldsymbol{0}$, $J \gets \emptyset$, $\boldsymbol{\ell}' \gets \boldsymbol{\ell}$, $\boldsymbol{u}' \gets \boldsymbol{u}$, $\boldsymbol{b}' \gets \boldsymbol{b}$\;
\While{$|J| < n$}{
Choose $i \in \bar{J} := \{1,\dots,n\} \setminus J$, $\bar{J} \gets \bar{J} \setminus \{i\}$\;
$\ell'_i \gets \displaystyle\max \left\{\ell_i', \max_{r:\,a_{ri} < 0} f(r, i, \bar{J}, \mv{\ell}', \mv{u}', \mv{b}') \right\}$\;
$u'_i \gets \displaystyle\min \left\{u_i', \min_{r:\,a_{ri} > 0} f(r, i, \bar{J}, \mv{\ell}', \mv{u}', \mv{b}') \right\}$\;
\If{$\ell'_i \leq u'_i$}{Sample $x^{k}_i \sim U\!\bigl(\ell'_i,\,u'_i\bigr)$}
\For{$r=1,\dots,m: a_{ri} \neq 0$}{
$b'_r \gets b'_r - a_{ri}\,x^{k}_i$\;
\For{$j \in \bar{J}$}{
$\ell'_j \gets \displaystyle\max \left\{\ell'_j, \max_{r:\,a_{rj} < 0} f(r, j, \bar{J}\setminus\{j\}, \mv{\ell}', \mv{u}', \mv{b}') \right\}$\;
$u'_j \gets \displaystyle\min \left\{u'_j, \min_{r:\,a_{rj} > 0} f(r, j, \bar{J}\setminus\{j\}, \mv{\ell}', \mv{u}', \mv{b}') \right\}$
}
}
$J \gets J \cup \{i\}$\;
}
$X \gets X \cup \{\boldsymbol{x}^{k}\}$\;
}
\end{algorithm}

Obtaining training labels for each sampled feasible first-stage investment decision presents a computational challenge. The expected production cost must be estimated using sample average approximation $\mathbb{E}_{\mv{\omega}} [Q(\mv{x}, \mv{\omega})] \approx \frac{1}{S} \sum_{s=1}^S Q(\mv{x}, \mv{\omega}^s)$ where each $Q(\mv{x}, \mv{\omega}^s)$ requires solving a production cost model over an operational horizon of length $|\mathcal{H}|$ (measured in chronological hours). This requires determining number of scenarios $S$ and operational horizon length $|\mathcal{H}|$. %
If they are chosen too small (e.g., $S=5$, $|\mathcal{H}|=12$), evaluating the production costs can be fast, but the estimates can be unreliable and in turn lead to poor investment decisions. On the contrary, if they are chosen too large (e.g., $S=100$, $|\mathcal{H}|=8760$), more accurate estimates are possible albeit at a steep computational price. Moreover, we can only know whether our parameters were sufficient in hindsight, and if not, we must regenerate all labels and retrain the surrogate from scratch.

To address these issues, we employ a hierarchical strategy for parameter selection motivated by recent empirical findings \cite{bylling2020impact}, which demonstrate that operational variability (temporal resolution) has a larger impact on reliability than scenario uncertainty (number of scenarios). Our algorithm proceeds in two nested loops:
\begin{itemize}
    \item Step 1 dynamically extends the operational horizon until the coefficient of variation of the estimated production cost falls below a given threshold~$\epsilon$.
    \item Step 2 increases the number of scenarios only if the resulting width of the $(1-\alpha)$ confidence interval for the expected cost $\widehat{\mathcal{Q}}$ exceeds some given tolerance~$\epsilon'$.
\end{itemize}

The proposed method, shown in Algorithm~\ref{alg:adaptive_labeling}, requires access to a scenario generator $\mathcal{P}(|\mathcal{H}|)$ that takes as input an operational horizon length $|\mathcal{H}|$ (measured in hours) and returns a randomly 
sampled scenario consistent with that horizon length, denoted as $\omega \sim \mathcal{P}(|\mathcal{H}|)$.

\begin{algorithm}[hbt!]
\small
\SetKwProg{Function}{function}{:}{end function}
\caption{Adaptive Two-Step Labeling Procedure}\label{alg:adaptive_labeling}
\DontPrintSemicolon
\KwIn{Sampled first-stage decision $\boldsymbol{x}$; initial number of scenarios $S_0 \ge 2$; initial operation horizon length {$H_0$}; horizon increment {$\Delta H$}; tolerances $\epsilon, \epsilon'$; confidence $\alpha$.}
\KwOut{Estimated expected production cost $\widehat{\mathcal Q}$.}
\Function{$\mathrm{Estimator} (\boldsymbol{x}, S_{\text{start}}, S_{\text{end}}, |\mathcal{H}|)$}{
$\text{Sum} \gets 0$; $\text{SumSq} \gets 0$\;
\For{$s \in \{S_{\text{start}}, \ldots, S_{\text{end}}\}$}{
Sample $\omega^s \sim \mathcal{P}(|\mathcal{H}|)$\;
$q_s \gets  Q(\boldsymbol{x},\boldsymbol{\omega}^{s})$\;
$\text{Sum} \gets \text{Sum}+q_s$; $\text{SumSq} \gets \text{SumSq}+q_s^2$\;
}
\Return $(\text{Sum}, \text{SumSq})$\;
}

$\widehat{\mathcal Q} \gets 0$, $S \gets S_0$, $|\mathcal{H}| \gets H_0$, $\delta \gets \infty$, $\bar{q} \gets \infty$,  $s^2 \gets \infty$\;
\While{$\delta > \epsilon$}{
$(\text{Sum}, \text{SumSq}) \gets {\mathrm{Estimator}(\boldsymbol{x}, 1, S, |\mathcal{H}|)}$\;
$\bar q \gets \text{Sum}/S$;\quad $s^2 \gets \frac{\text{SumSq} - \text{Sum}^2/S}{S-1}$\;
$\delta \gets \frac{\sqrt{s^2}}{\bar q}$\;
$|\mathcal{H}| \gets |\mathcal{H}|+\Delta H$\;
}
$r \gets z_{1-\alpha/2}\,\frac{\sqrt{s^2}}{\bar q \sqrt{S}}$\;
$\text{Sum}^* \gets \text{Sum}$;\quad $\text{SumSq}^* \gets \text{SumSq}$\;
\If{$r > \epsilon'$}{
$\widehat{S} \gets \max \!\left\{ \left\lceil \left(\frac{z_{1-\alpha/2}}{\epsilon'}\right)^2 \frac{s^2}{\bar q^2} \right\rceil,\, S \right\}$\;
\If{$\widehat{S} > S$}{
$(\text{Sum}', \text{SumSq}') \gets {\mathrm{Estimator}(\boldsymbol{x}, S{+}1, \widehat{S}, |\mathcal{H}|})$\;
$\text{Sum}^* \gets \text{Sum} + \text{Sum}'$\;
$\bar q \gets \text{Sum}^*/\widehat{S}$\;
}
}
$\widehat{\mathcal Q} \gets \widehat{\mathcal Q} + \bar q$\;
\Return $\widehat{\mathcal Q}$\;
\end{algorithm}

Our key hypothesis is that increasing the operational horizon length $|\mathcal{H}|$ reduces variance in production cost estimates by capturing temporal dependencies and averaging out short-term fluctuations in renewable generation and load. 
We provide evidence to substantiate this claim in Section~\ref{sec:results}.
Algorithm~\ref{alg:adaptive_labeling} therefore extends $|\mathcal{H}|$ until the coefficient of variation falls below some threshold $\epsilon$. In step 2, to decide a minimum sufficient number of scenarios $S$, we use the width of the confidence interval of the expected second-stage value function as a termination criterion. This is motivated by the standard techniques for evaluating solutions in stochastic programming \cite{shapiro2007tutorial}. Since the expected cost is estimated from a finite scenario sample, the confidence interval quantifies the precision of this estimate, and Step~2 augments the scenario count only when the current precision is insufficient.

As a general guideline, setting $\Delta H \le H_0$ avoids skipping horizon lengths that may already satisfy convergence, while a small $S_0 \ge 2$ suffices since the algorithm augments scenarios automatically. Setting $\epsilon \le \epsilon'$ is also advisable, as extending the horizon increases per-scenario cost whereas adding scenarios can be parallelized at fixed per-scenario cost.

For multi-horizon problems~\eqref{eq:mhsp_cep}, this algorithm can be applied independently to each investment period $t \in \mathcal{T}$, and the total estimated cost is obtained by summing across periods. This benefit leads to a significant computational advantage by allowing the label generations to be parallelized across both investment periods and scenarios. 

\subsection{Surrogate Model Architectures and Training}

We use the generated training data, $(\mv{x}^k, \widehat{\mathcal{Q}}^k)_{k=1}^K$, to learn a regression model.
We explore both simple linear models $\Phi^{\mathrm{LR}}$ and more complex neural network models $\Phi^{\mathrm{NN}}$.
Linear regression and variants (e.g., lasso) entail minimal computational effort for training and embedding (as we show below). However, they cannot capture the nonlinear nature of the expected value function. In contrast, feedforward neural networks with ReLU activation can potentially learn this complex nonlinear relationship, but typically require larger training datasets, computational resources and time. 
The choice of ReLU activations is important for the subsequent embedding of the trained model, as we discuss in the next subsection.
We explore the trade-off between model simplicity and predictive power in Section~\ref{sec:results}. 

\subsection{Surrogate Model Embedding}

The final step is to replace the expected value function $\mathcal{Q}(\mv{x})$ in equation~\eqref{eq:2sp_cep} with a trained regression model, $\Phi^{\mathrm{ML}}$. The resulting optimization problem is formulated as follows:
\begin{align}
    \mathop{\mathrm{minimize}} \;\quad & \mv{c}^{\top}_{\mathrm{inv}}\mv{x} + \mu  \label{eq:surrogate_obj} \\
    \text{subject to} \quad &\mv{A}\mv{x} \leq \mv{b}, \; \mv{\ell} \leq \mv{x} \leq \mv{u} \label{eq:feasible_region}\\
    & \mu \ge \Phi^{\text{ML}}(\boldsymbol{x}) \label{eq:surrogate_const} \\
    & \text{Constraints defining } \Phi^{\text{ML}} \label{eq:surrogate_ml} 
\end{align}

The trained model is embedded into the original problem as new constraints represented by equations~\eqref{eq:surrogate_const}--\eqref{eq:surrogate_ml}. These constraints involve only the epigraphical variable $\mu$ and the first-stage decisions $\mv{x}$, and do not modify the original feasible region defined by~\eqref{eq:feasible_region}. Embedding $\Phi^{\mathrm{LR}}$ in the optimization model requires only a single linear constraint, $\mu \geq \boldsymbol{\beta}^\top \boldsymbol{x} + \beta_0$, preserving the linear programming structure. On the contrary, embedding a feedforward ReLU network requires mixed-integer linear constraints; indeed, ReLU activations can be represented using binary variables and linear inequalities \cite{fischetti2018deep}.
For practical implementation, several open-source libraries facilitate the embedding of trained machine learning models into optimization problem, such as OMLT \cite{ceccon2022omlt}, PySCIPOPT\-ML \cite{turner2024pyscipoptmlembeddingtrainedmachine}, and Gurobi machine learning \cite{gurobiml}. %

We note that even though the original stochastic program is linear in the first and second-stage, solving it using conventional algorithms is not easy because of the large number of scenarios and operational hours. Therefore, even though we obtain a mixed-integer linear program (MILP) by embedding the trained feedforward ReLU model, the complexity of the MILP can be controlled by explicitly controlling the number of neurons or layers during neural network training.

\section{Computational Experiments} \label{sec:experiments}

\subsection{Benchmark Model}

We demonstrate AutoSCEP on the European Model for Power System Investment with Renewable Energy (EMPIRE), a large-scale capacity expansion planning model designed to determine capacity investments for the European power system from 2020 to 2060 \cite{backe2022empire}. The EMPIRE model has been widely used to inform energy investment strategies and conduct policy analyses \cite{durakovic2023powering, ahang2025investments}. The uncertain parameters in this model comprise hourly electricity demand profiles, generator availability (hourly capacity factors), and seasonal maximum production limits for
hydro generators across all nodes. These are jointly sampled from historical data to preserve spatial and temporal correlations.

The model is large: solving the deterministic EMPIRE model with a single scenario already requires two hours and 200GB of RAM, despite it being a linear program. To enable rapid experimentation and insights, we utilize two scaled-down versions for our analysis: intermediate- and reduced-size EMPIRE models, which we denote EMPIRE-med and EMPIRE-sml, respectively.

\begin{table}[hbt!]
\centering
\caption{Characteristics of the considered test systems}
\label{tab:problem_instance_compare}
\begin{tabular}{lccc}
\hline\textbf{}
\textbf{Parameter} &\textbf{EMPIRE-med}   & \textbf{EMPIRE-sml}\\
\hline
Periods  & 8 & 8 \\
Seasons  & 6 & 6 \\
Operational Hours & 336  & 240\\
Nodes  & 20 & 3 \\
Generators  & 373 & 65  \\
Transmission Lines  & 64 & 2 \\
Storage Units  & 38 & 5 \\
First-Stage Decisions  & 3,848 & 616 \\ 
Second-Stage Decisions & 1,056,384 & 130,560 \\
\hline
\end{tabular}
\end{table}

Table \ref{tab:problem_instance_compare} shows the sizes of the considered test systems. The EMPIRE-sml model has 616 decisions in the first stage and 130,560 in the second stage (per scenario). It corresponds to reductions of 84\% and 88\% compared to the EMPIRE-med model. %

\subsection{Comparative Baselines}

To evaluate AutoSCEP, we benchmark its performance against widely established algorithms for solving SP problems: Extensive Form (EF), Benders Decomposition (BD), and Progressive Hedging (PH). The EF is the deterministic equivalent of the stochastic program~\eqref{eq:2sp_cep}. Given a set of $S$ scenarios, the EF, which we denote as EF($S$), is formulated as follows:
\begin{subequations}
\label{eq:extensive_form}
\begin{align}
    \underset{\mv{x}, \{\mv{y}_s\}_{s=1}^S}{\textrm{minimize}} \quad & \mv{c}_{\textrm{inv}}^{\top}\mv{x} + \frac{1}{S}\sum_{s=1}^{S} \left( \sum_{h \in \mathcal{H}} c_{\textrm{prod}, \mv{\omega}^s}(\mv{y}_{hs}) \right) \\
    \text{subject to} \quad & \mv{A}\mv{x} \le \mv{b}, \quad \mv{\ell} \le \mv{x} \le \mv{u} \\
    & \mv{y}_s \in \mathcal{Y}(\mv{x}, \mv{\omega}^s), \quad \forall s \in \{1, \dots, S\}
\end{align}
\end{subequations}

Although the EF($S$) formulation is straightforward, it is often computationally intractable for large-scale problems. Decomposition algorithms like BD($S$) and PH($S$) overcome this challenge by creating subproblems for each scenario $s$, enabling them to handle a large number of scenarios efficiently. For the PH implementation, we adopted the parameter-free, adaptive penalty selection method \cite{watson2011progressive}. 

The convergence criteria was defined as the relative error difference between the non-anticipative solutions in each scenario and their across-scenario average solution, with a threshold of $10^{-3}$. We tested each baseline method on instances with $S \in \{5, 10, 20\}$ scenarios. For each case, we generated 10 distinct scenario sets using different random seeds. For the decomposition algorithm, we set time limits of 60, 300, 600, and 3600 seconds for EMPIRE-sml, and a time limit of six hours for EMPIRE-med.

\subsection{Implementation Details}

All optimization problems were solved using Gurobi 10.0.3 \cite{gurobi} with 1\% MIP gap tolerance. The experiments were conducted on Penn State's Roar HPC cluster, equipped with Intel Xeon CPU 6226R processors (180 CPUs) and 1TB of RAM. 
The baseline methods were implemented using mpisspy \cite{mpi-sppy}. Each baseline instance was allocated 64GB of RAM (400GB for EMPIRE-med) and multiple CPU cores per scenario for parallel decomposition. The data generation algorithm utilized 4GB RAM and 40 CPUs per training label. 

We used the following hyperparameter values for Algorithm~\ref{alg:adaptive_labeling}: coefficient of variation threshold $\epsilon=0.05$, confidence interval half-width tolerance $\epsilon'=0.1$, $S_0 = 5$, $H_0=6$, $\Delta H=6$, and confidence level $\alpha=0.05$.

For training Linear Regression (LR) and Multi-Layer Perceptron (MLP) with ReLU activations, we use scikit-learn 1.5.1 \cite{scikit-learn} and PyTorch 2.5.1 \cite{paszke2019pytorch}, respectively. 
The specific MLP architecture comprises two hidden layers with 16 and 8 neurons, respectively, with ReLU activations. Trained for up to 500 epochs with early stopping (patience of 20 epochs). We employed the Adam optimizer with a learning rate of 5e-5 and weight decay of 1e-5. Both models employ mean squared error loss function.

\section{Comparative Results}\label{sec:results}
We now examine three key questions:
\begin{itemize}
    \item How effective and efficient are our algorithms for training data generation, particularly in deciding the number of scenarios and operational hours?
    \item How effectively do the surrogate models balance computational efficiency with solution quality?
    \item How closely do the obtained first-stage capacity expansion decisions approximate the optimal solution?
\end{itemize}

All solutions were validated on 1000 independent out-of-sample scenarios that were not used during training. As an approximate measure of optimality, we calculate solutions by solving EF(100) for EMPIRE-sml and EF(20) for EMPIRE-med without any computational time restrictions.

\subsection{Adaptive Parameter Selection}

To evaluate the performance of Algorithm~\ref{alg:adaptive_labeling} in terms of adaptively finding an accurate number of scenarios and operational hours, we compare it with a fixed non-adaptive label generation approach with $S \in \{ 5, 10, 20, 30\}$ scenarios and $|\mathcal{H}| \in \{12, 24, 36, 48\}$ operational hours. %

\begin{figure*}[hbt!]
    \centering
    \includegraphics[width=1\linewidth]{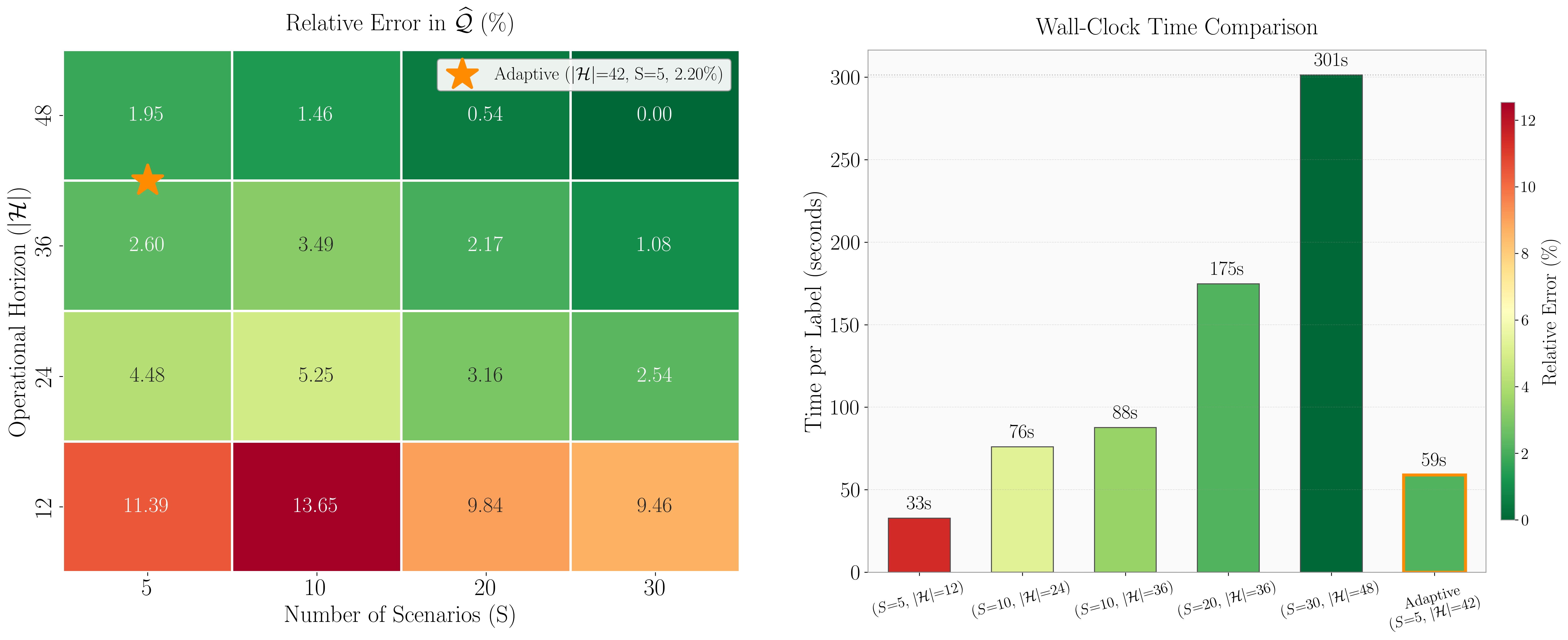}
    \caption{Comparison of fixed and adaptive parameter selection for generating a training label for a fixed $\mv{x}$. (Left)~Heatmap of relative error in $\widehat{\mathcal{Q}}$ with respect to the largest configuration ($S=30, |\mathcal{H}|=48$) for varying scenario counts $S$ and operational horizon lengths $|\mathcal{H}|$. The star marks the parameters selected by Algorithm~\ref{alg:adaptive_labeling} (2.2\% relative error). (Right)~Wall-clock time per label for representative fixed configurations and the adaptive method; bar color encodes relative error using the same colormap as the left panel.}
    
    \label{fig:param_comparison}
\end{figure*}

Figure~\ref{fig:param_comparison} compares the fixed-parameter and AutoSCEP in terms of accuracy (left) and efficiency (right). For this particular sampled decision, the heatmap shows that our approach converges to parameters that estimate the expected production cost $\widehat{\mathcal{Q}}(\mv{x})$ with high accuracy; i.e., it is almost equal to the cost obtained using the largest considered values of $S = 30$ and $|\mathcal{H}| = 48$.
At the same time, the right-hand side plot shows that when compared with the fixed parameter methods, our approach expends significantly less time in generating a near-identical label. These results demonstrate that our approach effectively balances between accuracy with computational tractability.

In the experiment above, Step~1 of Algorithm~\ref{alg:adaptive_labeling} (horizon extension) alone was sufficient to meet the precision threshold, so Step~2 (scenario augmentation) was not triggered. To verify that Step~2 activates when needed, we repeated the experiment with a smaller initial scenario count $S_0 = 2$. In this case, Algorithm~\ref{alg:adaptive_labeling} augments the scenarios to $S = 8$ and selects $|\mathcal{H}| = 36$, yielding an expected cost within 2.7\% of the largest configuration. This confirms that both steps of the algorithm contribute as needed to identify sufficient parameters.

\subsection{Capacity Expansion Solution Quality and Optimality}

We evaluate the performance of decomposition methods on both EMPIRE-sml and EMPIRE-med to understand the computational challenges that motivated our model reduction approach. BD fails to find feasible solutions, maintaining optimality gaps exceeding 18\% even after 3600 seconds. In contrast, PH converges rapidly on EMPIRE-sml, achieving optimality gaps below 1\% in 300 seconds across different scenario configurations. On the EMPIRE-med, PH also provides tight bounds with gaps around 4-5\% after 6 hours, though convergence is slower due to increased model complexity. We subsequently use PH as the validation benchmark for AutoSCEP.

We evaluate the surrogate model architectures and parameter selection strategies. Table~\ref{tab:model_performance} compares the performance of our proposed adaptive parameter selection algorithm against fixed parameter approaches for generating training datasets across different model architectures and dataset sizes.

\begin{table}[hbt!]
\centering
\caption{Comparison of adaptive and fixed parameter selection methods}
\label{tab:model_performance}
\resizebox{\columnwidth}{!}{%
\begin{tabular}{l l S[table-format=2.2] S[table-format=5] S[table-format=6] c}
\toprule
\textbf{Model} & \textbf{Method} & \textbf{MAPE (\%)} & \textbf{$L_T$ (s)} & \textbf{Total (s)} & \textbf{Speedup} \\
\midrule
\multirow{2}{*}{LR(1K)} & Adaptive & 10.79 & 7288 & 8512 & \multirow{2}{*}{4.31} \\
& Fixed & 10.79 & 35413 & 36697 & \\
\cmidrule(lr){1-6}
\multirow{2}{*}{MLP(1K)} & Adaptive & 8.18 & 7288 & 8536 & \multirow{2}{*}{4.30} \\
& Fixed & 9.47 & 35413 & 36717 & \\
\cmidrule(lr){1-6}
\multirow{2}{*}{LR(5K)} & Adaptive & 6.75 & 35150 & 41373 & \multirow{2}{*}{4.42} \\
& Fixed & 6.75 & 176820 & 183044 & \\
\cmidrule(lr){1-6}
\multirow{2}{*}{MLP(5K)} & Adaptive & 5.94 & 35150 & 41475 & \multirow{2}{*}{4.41} \\
& Fixed & 5.67 & 176820 & 183130 & \\
\bottomrule
\end{tabular}
}
\begin{flushleft}
\footnotesize
\textbf{Note:} Model designations such as MLP(1K) specify the model type and the dataset size (e.g., 1,000 samples). ``Total'' time is the sum of sampling, labeling ($L_T$), and model training times. "Speedup" is the ratio of the total execution time of the Fixed method to the Adaptive method. All metrics are averaged over 10 different sample datasets.
\end{flushleft}
\end{table}

Table~\ref{tab:model_performance} demonstrates the effectiveness of our adaptive parameter selection algorithm. The results from both parameter selection methodologies demonstrate that predictive accuracy improves with larger dataset sizes. Across all configurations, MLP models consistently outperform the LR models. With the 5K-sample dataset, the MLP model achieves a low mean absolute percentage error (MAPE) of 5.67\% with adaptive parameter selection method compared to 5.94\% with fixed method. The most significant finding is the computational efficiency of our proposed algorithm: it achieves comparable accuracy to the fixed parameter methods while reducing total dataset generation time by approximately 77\%.

We now validate the solutions generated by AutoSCEP against baseline algorithms. Table~\ref{tab:comprehensive_results} presents comprehensive comparison of solution quality and computational performance across baseline methods and AutoSCEP.

\begin{table}[hbt!]
\centering
\caption{Performance trade-offs between solution quality and time across different methods.}
\label{tab:comprehensive_results}
\begin{tabular}{ll
                S[table-format=2.2]
                S[table-format=1.2]
                S[table-format=5.2]}
\toprule
& & \multicolumn{2}{c}{\textbf{\makecell{Solution \\ Quality}}} & {\textbf{\makecell{Solution \\ Time}}} \\
\cmidrule(lr){3-4} \cmidrule(lr){5-5}
\textbf{Model} & \textbf{Method} & {\textbf{Gap (\%)}} & {\textbf{CV (\%)}} & {\textbf{Time (s)}} \\
\midrule
\multirow{14}{*}{\textbf{EMPIRE-sml}}
 & EF(5)   & 2.24 & 1.14 & 722.67 \\
 & EF(10)  & 0.91 & 0.34 & 1387.26 \\
 & EF(20)  & 0.52 & 0.18 & 3100.26 \\
 \cmidrule(lr){2-5} 
 & PH(5)   & 1.88 & 0.97 & 300.00 \\
 & PH(10)  & 1.21 & 0.44 & 300.00 \\
 & PH(20)  & 0.58 & 0.11 & 300.00 \\
 \cmidrule(lr){2-5}
 & F-LR(1K)  & 5.48 & 2.91 & 0.86 \\
 & F-MLP(1K) & 2.79 & 1.43 & 0.97 \\
 & F-LR(5K)  & 1.39 & 1.21 & 1.15 \\
 & F-MLP(5K) & 1.90 & 1.04 & 1.09 \\
 \cmidrule(lr){2-5}
 & A-LR(1K)  & 5.48 & 2.91 & 1.19 \\
 & A-MLP(1K) & 2.91 & 1.38 & 1.22 \\
 & A-LR(5K)  & 1.96 & 1.11 & 1.78 \\
 & A-MLP(5K) & 1.61 & 0.73 & 1.51 \\
\midrule
\multirow{4}{*}{\textbf{EMPIRE-med}}
 & PH(10)  & 11.38 & 2.26 & 21600.00 \\
 & PH(20)  & 11.06 & 2.51 & 21600.00 \\
 \cmidrule(lr){2-5}
 & A-LR(500) & 11.28 & 5.03 & 5.81 \\
 & A-MLP(500)  & 8.14 & 2.44 & 54.57 \\
\bottomrule
\end{tabular}
\begin{flushleft}
\footnotesize
\textbf{Note:} Gap = Optimality Gap calculated against EF(100) for EMPIRE-sml and EF(20) for EMPIRE-med. CV = Coefficient of Variation. Methods prefixed with ``A'' denotes AutoSCEP and ``F'' denotes Neur2SP with a specified model trained using datasets generated with fixed parameters.
\end{flushleft}
\end{table}

As shown in Table \ref{tab:comprehensive_results}, AutoSCEP achieves dramatic computational speedups compared to traditional algorithms. On EMPIRE-sml, the A-MLP(5K) method yields a solution with 1.6\% optimality gap in less than 2~seconds. %
Critically, AutoSCEP also outperforms its fixed parameter counterparts: A-MLP(5K) yields a lower optimality gap and variance compared to F-MLP(5K). Moreover, the trade-off between model simplicity and accuracy is evident. The LR models (A-LR and F-LR) achieve faster solution times but with higher optimality gaps and CV, while MLP models (A-MLP and F-MLP) provide lower gaps and CV with slightly increased computational time across all experimental configurations. This demonstrates that embedding a more complex neural network provides better solutions within practical runtimes.

The fixed-parameter methods (F-LR, F-MLP) are omitted from the EMPIRE-med analysis. As shown in Table~\ref{tab:model_performance}, the fixed-parameter approach is roughly four times slower than the adaptive method for generating training labels on EMPIRE-sml. Because EMPIRE-med has approximately eight times more second-stage decision variables per scenario, the estimated wall-clock time for generating even 500 fixed-parameter labels exceeds the available computational budget, further motivating the adaptive approach. The scalability of AutoSCEP becomes particularly evident on the EMPIRE-med model. Whereas PH requires 6~hours to achieve 11\% gap, AutoSCEP with MLP (A-MLP(500)) finds a better solution with 8\% gap and comparable CV of 2.4\%.

\subsection{Structural Reliability of Investment Plans}

\begin{figure*}[hbt!]
    \centering
    \includegraphics[width=1\textwidth]{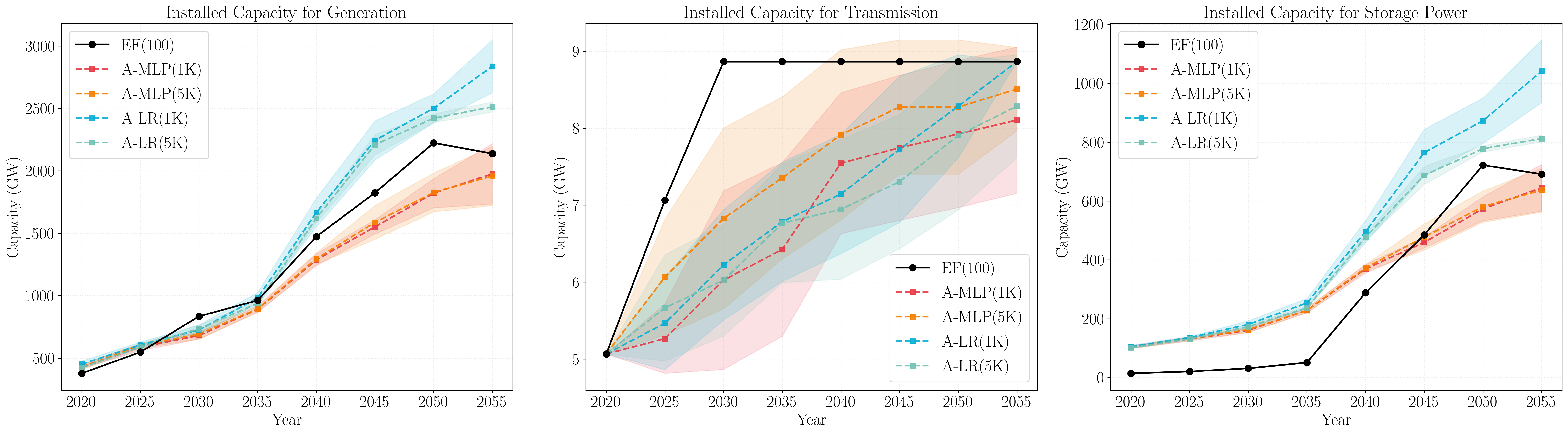}
    \caption{Cumulative installed capacity by technology for the EMPIRE-sml case study; shaded regions represent 95\% confidence intervals across ten solutions.}
    \label{fig:structural_comparison}
\end{figure*}

To evaluate our solutions in terms of their structure similarity to the optimal solution, Figure~\ref{fig:structural_comparison} plots the temporal capacity expansion decisions suggested by the different methods across four technology types. The optimal strategy appears to involve: (i) a rapid investment in transmission capacity during the early planning period; (ii) a steady investment in generation capacity; and (iii) delaying investment in energy storage systems until after 2035.

We find that the A-MLP(5K) method successfully learns the optimal patterns, particularly in generation and the timing of storage capacity investments. Although it does not perfectly mimic the transmission investment strategy, the A-MLP(5K) method still outperforms the other models, and its overall strategy comes quite close to the EF(100) benchmark.

In contrast, the other models do not mimic the optimal strategy of EF(100). In particular, the A-LR(5K) model suggests a structurally different plan: delaying the early period investment in transmission, despite over-investing in generation and storage.

\section{Conclusion} \label{sec:conclusions}

This paper presented AutoSCEP, a practical framework that lowers the computational barrier to large-scale capacity expansion planning. AutoSCEP replaces ad hoc choices of scenario count and operational horizon with an adaptive, statistically grounded procedure that selects the minimum sufficient values to meet a prescribed precision. AutoSCEP uses these labels to train linear and neural surrogates of expected production cost and embeds them within the planning model. The surrogates can be reused across candidate plans, which limits repeated second-stage simulation and reduces tuning and implementation effort.
On the continental-scale EMPIRE system, AutoSCEP attained less than 2 percent optimality gap on a reduced model and about 8 percent on a large model. Under equal wall-clock budgets that include data generation, training, and solve times, it remains competitive with parallel progressive hedging. The speed of the workflow enables rapid evaluation of alternative investment strategies in response to policy changes, while preserving high-resolution uncertainty modeling at realistic system scales.

Our approach has several limitations that suggest directions for future work. First, the trained surrogate is specific to the second-stage operational model used during data generation; if the pipe-and-bubble network model in EMPIRE were replaced with, say, a nodal unit commitment formulation, or if the network topology or resolution changed, the surrogate would need to be retrained. 
The proposed algorithms, however, carry over without modification to all such cases, and exploring more detailed CEP models is a natural next step. Second, the hierarchical strategy in Algorithm~\ref{alg:adaptive_labeling} extends the operational horizon before augmenting the scenario count, an ordering motivated by~\cite{bylling2020impact} for a single-node system. Our results in Figure~\ref{fig:param_comparison} support this ordering within the multi-node EMPIRE model, but it may not hold universally, and further investigation on diverse test systems is warranted. Third, the decomposition benchmarks in Section~\ref{sec:results} use the default initialization and settings of the mpi-sppy framework~\cite{mpi-sppy}, to which their performance is known to be sensitive~\cite{jacobson2024computationally, watson2011progressive}; better tuning could improve their results. AutoSCEP avoids such initialization sensitivity but introduces its own hyperparameters (e.g., the number of ReLU neurons and layers) that govern a tradeoff between surrogate accuracy and solve time. 
A systematic study of this tradeoff along with graph-based surrogates for improved scalability are promising avenues for future research.

\bibliographystyle{IEEEtran}
\bibliography{citation.bib}

\end{document}